\documentclass[11pt]{amsart}
\usepackage{amsmath,amssymb, amsthm, amsfonts}
\usepackage{verbatim}
\let\d=\partial

\let\wt=\widetilde

\def\cC{{\mathcal C}}

\def\cO{{\mathcal O}}
\def\cP{{\mathcal P}}
\def\cQ{{\mathcal Q}}
\def\cR{{\mathcal R}}

\def\N{{\mathbb N}}
\def\R{{\mathbb R}}

\def\bv{{\bar v}}

\def\virgp{\raise 2pt\hbox{,}}
\def\cdotpv{\raise 2pt\hbox{;}}
\def\Id{\mathop{\rm Id}\nolimits}

\def\adj{\mathop{\rm adj}\nolimits}

\def\div{ \hbox{\rm div}\,  }
\def\divv{ \hbox{\rm div}_{\!v}\,  }

\newcommand{\Int}{\displaystyle \int}

\newcommand{\ds}{\displaystyle}

\newcommand{\with}{\quad\hbox{with}\quad}
\newcommand{\andf}{\quad\hbox{and}\quad}

\def\dv{\delta\!v}

\def\deta{\delta\!\eta}

\newtheorem{thm}{Theorem}[section]

\newtheorem{rmk}{Remark}[section]

\newtheorem{prop}{Proposition}[section]

\newcommand{\ben}{\begin{eqnarray}}
\newcommand{\een}{\end{eqnarray}}
\newcommand{\beno}{\begin{eqnarray*}}
\newcommand{\eeno}{\end{eqnarray*}}
\numberwithin{equation}{section}

\begin{document}
\title[ Lorentz spaces and pressureless gases]
{Lorentz spaces in action on pressureless systems arising from  models of collective behavior}
\author[R. Danchin]{Rapha\"{e}l Danchin}
\address[R. Danchin]{Univ Paris Est Creteil, CNRS, LAMA, F-94010 Creteil, France,
Univ Gustave Eiffel, LAMA, F-77447 Marne-la-Vall\'ee, France,
Sorbonne Universit\'e and Universit\'e de Paris, CNRS, Inria, Laboratoire Jacques-Louis Lions (LJLL), F-75005 Paris, France} \email{raphael.danchin@u-pec.fr}
\author[P.B. Mucha]{Piotr Bogus\l aw Mucha}
\address[P.B. Mucha]{Instytut Matematyki Stosowanej i Mechaniki, Uniwersytet Warszawski, 
ul. Banacha 2,  02-097 Warszawa, Poland} 
\email{p.mucha@mimuw.edu.pl}
\author[P. Tolksdorf]{Patrick Tolksdorf}
\address[P. Tolksdorf]{Institut f\"ur Mathematik, Johannes Gutenberg-Universit\"at Mainz, Staudingerweg 9, 55099 Mainz, Germany}
\email{tolksdorf@uni-mainz.de}
\dedicatory{Zu Ehren von Matthias Hieber}
\begin{abstract}
We are concerned with global-in-time existence and uniqueness results for 
models of pressureless gases that come up 
in  the description of phenomena in astrophysics or collective behavior. 
 The initial data are rough: in particular, the density  is only bounded. 
Our results are based on interpolation and parabolic maximal regularity, 
where Lorentz spaces play a key role.  
We establish a novel maximal regularity estimate for parabolic systems in
$L_{q,r}(0,T;L_p(\Omega))$ spaces.

 \end{abstract}
\maketitle




\section{Introduction}
 We  are concerned with   models coming from a special type of hydrodynamical systems, that  do not include the effects of the internal pressure. 
The simplest example  is the motion of dust, that is, of free particles 
evolving in the space like e.g.  in astrophysics~\cite{DBD}, or in multi-fluid systems~\cite{Ber,BChZ}.  Leaving the world of inanimate matter,
 one can also mention models that describe  collective behavior, 
 where particles or rather agents exhibit some intelligence, and 
 for which having a force  like  internal pressure is not so natural. 
 A well-known example in this area is given by  the equations of traffic flow~\cite{BDD,HS}, where particles are just cars.
 \smallbreak
  In order to specify and understand this class of models, let us  go back to the kinetic description
of a collective behavior.  
  Consider  equations of the following form
 \begin{equation}\label{x1}
  f_t +v\cdot \nabla_x f + \div_{\!v} K(f)f=0 \mbox{ \ \ in \ } (0,T)\times \R^d_x \times \R^d_v,
 \end{equation}
 where $f=f(t,x,v)$ is a distribution function of a gas in the phase-space. 
 Classically, if the operator~$K$ comes from the Poisson potential, then we find the Vlasov system.   If  taking a less singular  operator, then one may obtain for example the Cucker-Smale system 
 that  models collective behavior like flocking of birds~\cite{CS}.
  \smallbreak
Assuming a very special form of $f$, the so-called
 \emph{mono-kinetic ansatz}, one can  pass formally 
 from the kinetic model  \eqref{x1} to the hydrodynamical system, putting just 
 \begin{equation}\label{x2}
  f(t,x,v)=\rho(t,x) \delta_{v=u(t,x)}.
 \end{equation}
 This  amounts to saying that the   distribution of the gas
 under consideration is located on  the curve $v=u(t,x)$.
Although one cannot expect this simplification to be a correct description 
of a gas, it is  relevant for modelling  some collective behavior phenomena as  
one can expect a crowd of individuals   to have the same speed (or tendency)
at  one point~\cite{MT,OPA}.
\smallbreak 
 Formally, plugging~\eqref{x2} in~\eqref{x1} leads to the following general form of hydrodynamical system: 
\begin{equation}\label{x3}
 \begin{array}{l}
  \rho_t + \div(\rho u)=0,\\[1ex]
  \rho u_t + \rho u\cdot\nabla u = A(\rho,u).
 \end{array}
\end{equation}
The two  equations may be seen  as  the mass and momentum balances, respectively. 
If $A\equiv 0,$ then one just recovers  the pressureless compressible Euler system, and there is  no interaction whatsoever between the individuals. 
Relevant examples where $A\not\equiv0$ can be found in~\cite{MM,ST,TT}. 
\medbreak
In our note, we would like to put our attention on the following two cases:
\begin{equation}\label{x4}
 A(\rho,u)= \mu \Delta u \ \hbox{ or }\ A(\rho,u)=\mu \Delta u + \mu'\nabla \div u,\qquad
 \mu>0,\quad \mu+\mu'>0.
\end{equation}
The first case is  a viscous regularization of~\eqref{x3} that can be viewed as a simplification of the Euler alignment system. It corresponds to the   hydrodynamical version of the Cucker-Smale model, namely
\begin{equation}\label{x5}
 \begin{array}{l}
  \rho_t +\div(\rho u)=0,\\
 \displaystyle \rho u_t + \rho u \cdot \nabla u = \int_{\R^d} \frac{u(t,y)-u(t,x)}{|x-y|^{d+\alpha}}
  \rho(t,y)\rho(t,x) dy,\qquad\qquad\alpha\in(0,2).
 \end{array}
\end{equation}
 The right-hand side of~\eqref{x5} involves the fractional Laplacian 
$(-\Delta)^{\alpha/2}$ (see details in~\cite{DMPW,Kie})
and the first case in~\eqref{x4} thus meets $\alpha = 2$.
 The second case of~\eqref{x4} is  the Lam\'e operator that can be obtained from 
 the Vlasov-Boltzmann equation (for more explanation, one may refer to the introduction of~\cite{Laura}).
\smallbreak
 The form of~\eqref{x4}  does not take into account  the effects of internal pressure. 
 From the mathematical viewpoint, the lack of the pressure term $P$ causes serious problems. 
 In particular, all techniques for the compressible viscous systems
 based on the properties of the so-called effective viscous flux, 
 namely $F:=  \div u - P,$ 
which has better regularity than $\div u$ and $P$ taken separately, are bound to fail.
Recall that  using $F$ is one of the keys  to the theory of weak solutions of the compressible Navier-Stokes equations~\cite{Hoff1,Hoff2,Lions,Fei,FM},  as it allows to exhibit compactness properties of the set of weak solutions. In the theory of regular solutions~\cite{MN,D00,Mu03},  the effective viscous flux  provides the decay properties for the density
that are needed for establishing global existence for small data. 

In the case of pressureless systems,  there is no such a possibility,  so that we need to resort to  more sophisticated techniques to control the density. 
This may partially explain the reason why  the mathematical theory of 
pressureless models is poorer than the classical one.
\smallbreak 
 The aim of this note is to present a novel
  technique coming from the maximal regularity theory  for analytic semi-groups,
   allowing  to prove global-in-time properties of solutions to~\eqref{x3},~\eqref{x4}.
  It will enable us to show 
  existence and uniqueness results under 
  rough assumptions on the density (only bounded), \emph{even though
  one cannot take advantage of the effective viscous flux.}
More precisely, by combining interpolation arguments, subtle
 embeddings,  suitable time weighted norms and
 the magic properties of Lorentz spaces, 
 we succeed in obtaining the $L_1(\R_+;L_\infty)$ regularity for  the gradient 
 of the solution to the linearized momentum equation in~\eqref{x3} 
 and, eventually,  produce global-in-time strong solutions. Our main results concern global-in-time solvability for the two dimensional case for large velocity, and the three dimensional case
  in the small data regime.


 \section{Functional framework} 

This note aims at presenting an  interesting application of Lorentz spaces for 
parabolic type systems. 
Lorentz spaces can be defined on any measure space $(X,\mu)$
 via  real interpolation between the classical Lebesgue spaces, as follows:
\begin{equation} L_{p,r}(X,\mu) := (L_\infty(X,\mu),L_1(X,\mu))_{1/p,r}
\ \hbox{ for }\ p\in (1,\infty)\andf r \in [1,\infty].
\end{equation}
 Lorentz spaces may be endowed with the following (quasi)-norm (see, e.g.,~\cite[Prop. 1.4.9]{grafakos}):
\begin{equation}
\|f\|_{L_{p,r}}:=\left\{\begin{array}{ll}
p^{\frac1r}\biggl(\Int_0^\infty  \bigl(s\, |\{|f|>s\}|^{\frac1p}\bigr)^r\,\frac{ds}s\biggr)^{\frac1r}
&\hbox{if }\ r<\infty,\\
\ds\sup_{s>0} s\, |\{|f|>s\}|^{\frac1p}&\hbox{if }\ r=\infty.\end{array}\right.
\end{equation}
The reason for the pre-factor $p^{\frac1r}$ is to have $\|f\|_{L_{p,p}}=\|f\|_{L_p},$
according to Cavalieri's principle.
%
%
The following classical properties of Lorentz spaces 
may be found in, e.g.,~\cite{Be,grafakos}:
\begin{itemize}
\item[--] Embedding : $L_{p,r_1}\hookrightarrow L_{p,r_2}$ if $r_1\leq r_2,$ 
and $L_{p,p}=L_p.$ 
\item[--] H\"older inequality :  for $1<p,p_1,p_2<\infty$ and
$1\leq r,r_1,r_2\leq\infty$, we have
\begin{equation}\label{eq:holder}
\|fg\|_{L_{p,r}}\lesssim \|f\|_{L_{p_1,r_1}} \|g\|_{L_{p_2,r_2}}\with
\frac 1p=\frac1{p_1}+\frac1{p_2}\andf 
\frac 1r=\frac1{r_1}+\frac1{r_2}\cdotp
\end{equation}
Inequality \eqref{eq:holder}  still holds for couples $(1,1)$ and $(\infty,\infty)$ with the convention
$L_{1,1}=L_1$ and $L_{\infty,\infty}=L_\infty.$
\item[--]  For any $\alpha>0$ and nonnegative measurable function $f,$ we have 
$$\|f^\alpha\|_{L_{p,r}}=\|f\|_{L_{p\alpha,r\alpha}}^\alpha.$$
\end{itemize}
It must be highlighted that Lorentz spaces have the following two interesting
and useful properties:
\begin{itemize}
 \item [--] If $f\in L_\infty(\Omega)$ and $g\in L_{p,q}(\Omega)$, then
 $fg \in L_{p,q}(\Omega)$ for all $p,q\in [1,\infty]$.
 \item [--] If $\Omega \subset \R^d$ is open and if $f\in L_{d,1}(\Omega)$ and $\nabla f \in L_{d,1}(\Omega)$,
 then $f \in C_b(\Omega)$.
\end{itemize}
However the class of Lorentz spaces possesses also ``false friends":
\begin{itemize}
 \item [--] There is no ``Fubini property",  that is  if  $(I,J)$ is a couple of nontrivial intervals, then 
 (see~\cite{Cwickel}):
 $$  L_{p,q}(I\times J)\neq L_{p,q}(I;L_{p,q}(J)) \quad\hbox{whenever }\ p \neq q.$$ 
 \item [--]  Being not reflexive, the space $L_{p,1}(\Omega)$ does not have the UMD property, 
  and   the general theory of deducing maximal regularity in $L_q(0,T;L_{p,1}(\Omega))$ via $\cR$-boundedness and Fourier multiplier theory developed in~\cite{Denk} cannot be applied.
\end{itemize}

Owing to their  definition  by real interpolation, 
Lorentz spaces have some connections with  homogeneous and nonhomogeneous Besov spaces. 
Recall that those spaces  may be defined respectively (for all $s\in(0,1)$ and $1\leq p,q\leq\infty$) by 
\begin{equation}
 \dot B^s_{p,q}(\Omega)=(L_p(\Omega);\dot W^1_p(\Omega))_{s,q}\andf
 B^s_{p,q}(\Omega)=(L_p(\Omega); W^1_p(\Omega))_{s,q}.
\end{equation}
For more properties of Besov spaces, we refer to \cite{BCD}.
\medbreak
The main results of the paper strongly rely on  a family of maximal regularity estimates
for the heat equation, where the time regularity is measured in Lorentz spaces
of type $L_{q,r}(0,T;L_p)$. Those  estimates come up as  a consequence of $L_q(0,T;L_p)$ maximal regularity estimates from the general theory \cite{Amann,Denk,Giga,Pia}.
\begin{prop}\label{p:MR}  Let $1<p,q<\infty$ and $1\leq r \leq \infty.$ Then, for any $u_0\in\dot B^{2-2/q}_{p,r}(\R^d)$ and $f\in L_{q,r}(\R_+;L_p(\R^d)),$
the following heat equation: 
\begin{equation}\label{heat}
\begin{aligned}
u_t-\mu\Delta u&=f\quad&\hbox{in }\ \R_+\times\R^d,\\
u|_{t=0}&=u_0&\quad\hbox{in }\ \R^d\end{aligned}
\end{equation}
has a unique solution in the space\footnote{Only weak continuity if $r=\infty$.}  
$$
\dot W^{2,1}_{p,(q,r)}(\R^d\times\R_+):=\Bigl\{u\in \cC_b\bigl(\R_+;\dot B^{2-2/q}_{p,r}(\R^d)\bigr):
u_t,\nabla^2 u\in  L_{q,r}(\R_+;L_p(\R^d))\Bigr\},$$
and the following inequality holds true:
\begin{multline}\label{eq:maxreg1}
\mu^{1-1/q}\|u\|_{L_\infty(\R_+;\dot B^{2-2/q}_{p,r}(\R^d))} + 
 \|u_t,\mu \nabla^2 u\|_{L_{q,r}(\R_+;L_{p}(\R^d))} \\ 
 \leq C\left(\mu^{1-1/q}\|u_0\|_{\dot B^{2-2/q}_{p,r}(\R^d)} + \|f\|_{L_{q,r}(\R_+;L_{p}(\R^d))}\right)\cdotp
\end{multline}
Furthermore, if $2/q + d/p>2$, then for all $q < s<\infty$ and $p\leq m$ such that $1+\frac d2(\frac1m-\frac1p)>0,$  interrelated by 
$$\frac d{2m}+\frac1s=\frac1q+\frac d{2p}-1,$$ it holds that
\begin{align*}
 \mu^{1+\frac1s-\frac1q}\|u\|_{L_{s,r}(\R_+;L_m(\R^d))} \leq C \Big( \mu^{1-1/q}\|u\|_{L_\infty(\R_+;\dot B^{2-2/q}_{p,r}(\R^d))} + \|u_t,\mu \nabla^2 u\|_{L_{q,r}(\R_+;L_{p}(\R^d))} \Big)\cdotp
\end{align*}
\end{prop}

The proof of the above result is given in the appendix.

\section{Results} 

Let us first  present our results pertaining to  the case where $A(\rho,u)=\mu \Delta u + \mu'\nabla \div u$
if the gas domain  is the whole plane $\R^2.$ So, we consider: 
 \begin{equation}\label{eq:presless}
  \begin{array}{lr}
   \rho_t +\div(\rho u)=0 & \mbox{ in }\ \R_+\times \R^2,\\[1ex]
   \rho u_t +\rho u \cdot \nabla u = \mu \Delta u + \mu'\nabla \div u
   & \mbox{ in }\ \R_+\times \R^2,\\[1ex]
   \rho|_{t=0} =\rho_0,\qquad u|_{t=0}=u_0 & \mbox{ at }\ \R^2.
  \end{array}
 \end{equation}
Following recent results of the first two authors in~\cite{DM-adv,DM-TJM,DM-Largenu}
(in different contexts, though), we strive for global results for 
general initial velocities  provided the volume (bulk) viscosity $\mu'$
is large enough. As our approach is based on a perturbative method, 
we need moreover the density to be close to some positive constant.
\medbreak
 Our solution space  will be the set $\dot W^{2,1}_{4/3,(4/3,1)}(\R^2\times\R_+)$ of functions $z:\R_+\times\R^2\to\R$
such that  
$$z\in\cC_b(\R_+;\dot B^{1/2}_{4/3,1})\andf
\partial_tz,\, \nabla^2_xz \in L_{4/3,1}(\R_+; L_{4/3}(\R^2)).$$
As shown in the appendix,   the corresponding trace space on fixed times is 
the homogeneous Besov space $\dot B^{1/2}_{4/3,1}(\R^2),$ 
which is `critical' in terms of regularity, and 
 embedded in  $L_2(\R^2).$ 
 \medbreak
 Our first statement  is a  global existence and uniqueness result for~\eqref{eq:presless} for large data in the two dimensional case, provided the volume viscosity $\nu$ is sufficiently large.
\begin{thm}\label{thm:pressureless}
Let  $M\geq1.$  Denote by $\cP$ and $\cQ$ the Helmholtz projectors
 on divergence-free and potential vector fields, and set $\nu:=\mu+\mu'.$

There exist two constants $c$ and $C$ (independent of $M$) such that if  the initial density $\rho_0$ satisfies
 \begin{equation}\label{eq:smallrho0} 
 \|\rho_0-1\|_{L_\infty(\R^2)}\leq c
 \end{equation}
and if the initial velocity  $u_0$ has components in $\dot B^{1/2}_{4/3,1}(\R^2)$ 
and satisfies,  \begin{equation}\label{eq:data}
 C\bigl(\|\cP u_0\|_{\dot B^{1/2}_{4/3,1}(\R^2)}\!+\!(\nu/\mu)^{1/4}\|\cQ u_0\|_{\dot B^{1/2}_{4/3,1}(\R^2)}\bigr)\exp\bigl\{C\mu^{-1}\|u_0\|_{L^2(\R^2)}^2\bigr\}\leq M\mu
\end{equation}
then, for all $\nu\geq \nu_0:= \max(C M^{1 / 2} e^{CM^4} , \mu),$
  System~\eqref{eq:presless} admits a global finite energy solution $(\rho,u)$
 with  $\rho\in\cC_w^*(\R_+;L_\infty(\R^2))$ and  $u\in \dot W^{2,1}_{4/3,(4/3,1)},$
 satisfying
  \begin{equation}\label{eq:smallrho1}  \|\rho-1\|_{L_\infty(\R_+\times\R^2)}\leq 2c
  \andf \|u\|_{\dot W^{2,1}_{4/3,(4/3,1)}(\R^2\times\R_+)}\leq C\|u_0\|_{\dot B^{1/2}_{4/3,1}(\R^2)}.\end{equation}
 Furthermore,  $\nabla u\in L_1(\R_+;L_\infty(\R^2))$  and the following decay property holds:
 \begin{multline}\label{eq:decay}
 \|t\cP u\|_{L_\infty(\R_+;\dot B^{3/2}_{4,1}(\R^2))} +\Bigl(\frac{\nu}{\mu}\Bigr)^{3/4}
\|t\cQ u\|_{L_\infty(\R_+;\dot B^{3/2}_{4,1}(\R^2))} \\+
\|(tu)_t,\nabla^2 tu\|_{L_{4,1}(\R_+;L_{4}(\R^2))} 
+\frac{\nu}{\mu}\|\nabla \div tu
\|_{L_{4,1}(\R_+;L_4(\R^2))} \leq  e^{CM^4}\cdotp\end{multline}
 Finally,   the solution $(\rho,u)$ is unique in the above regularity class.
 \end{thm}
 Some comments are in order:
 \begin{itemize}
 \item[--] Condition~\eqref{eq:data} means that global existence holds true 
 for large $\nu,$ provided the divergence part of the velocity is  $\cO(\nu^{-1/4}).$
 A similar restriction (with other exponents, though), was found in 
 our prior works dedicated to the global existence of strong solutions for
 the compressible Navier-Stokes equations with increasing pressure law~\cite{DM-adv,DM-TJM,DM-Largenu}.
 \item[--]  The above statement involves 
 only quantities that are scaling invariant for System~\eqref{eq:presless}.
 \item[--]  Having  Lorentz spaces with 
 second index equal to $1$  allows us to capture  limit cases of Sobolev embeddings. 
  Other choices than  $L_{4/3,1}$ and $L_{4,1}$   might be possible. 
 \end{itemize}
 \medbreak
  In the three-dimensional case,  the energy space $L_2(\R^3)$ is super-critical by 
 half a derivative, and   there is no chance (so far)
  to prove a general result for  large data, assuming only that one of the viscosity coefficients is large. 
  Therefore,  for simplicity, we  focus on the first  case of~\eqref{x4},  that is on 
   the following system: 
   \begin{equation}\label{eq:NS3D}
  \begin{array}{lr}
   \rho_t+\div (\rho u)=0 & \hbox{ in }\ \R_+\times\R^3,\\[1ex]
   \rho u_t +\rho u \cdot \nabla u = \mu \Delta u & \hbox{ in }\ \R_+\times\R^3,\\[1ex]
   u|_{t=0}=u_0, \qquad \rho|_{t=0}=\rho_0 & \hbox{ at }\ \R^3.
  \end{array}
 \end{equation}
 To simplify our analysis, we  choose a functional 
 framework for the velocity that is well beyond critical regularity (actually, 
 we ask for one more derivative).   In the same spirit as in the 2D case, 
 we will use the following class of functional spaces:
 $$\dot W^{2,1}_{p,(q,1)}(\R^3\times\R_+):=\Bigl\{u\in \cC_b\bigl(\R_+;\dot B^{2-2/q}_{p,1}(\R^3)\bigr):
u_t,\nabla^2 u\in  L_{q,1}(\R_+;L_p(\R^3))\Bigr\}\cdotp$$
 Our global existence result in the three dimensional case reads:
 \begin{thm}\label{thm:NS3D}
 Take initial data  $\rho_0 \in L_\infty(\R^3)$ and  $u_0\in \dot B^{3/5}_{10/7,1}(\R^3)\cap\dot B^{6/5}_{5/2,1}(\R^3).$ 
There exists a constant $c>0$ such that if
\begin{equation}\label{eq:smallness}
\|\rho_0-1\|_{L_\infty(\R^3)} <c \andf
\|u_0\|_{\dot B^{6/5}_{5/2,1}(\R^3)}^{1/3} \|u_0\|_{\dot B^{3/5}_{10/7,1}(\R^3)}^{2/3}<c\mu,
\end{equation}
then~\eqref{eq:NS3D} has a global-in-time unique solution $(\rho,u)$  with 
 $$\rho\in\cC_w(\R_+;L_\infty(\R^3))\andf u \in (\dot W^{2,1}_{5/2,(5/2,1)}\cap  \dot W^{2,1}_{10/7,(10/7,1)})(\R^3\times\R_+).$$
Furthermore, that solution has finite energy, 
\begin{equation}\label{eq:smallrho}
\|\rho-1\|_{L_\infty(\R^3)} < 2c,
\end{equation}
 functions $tu$  and $\nabla u$  belong  to $\dot W^{2,1}_{10/3,(10/3,1)}(\R^3\times\R_+)$ and  $L_1(\R_+;L_\infty(\R^3)),$ respectively, and 
 the following inequalities are valid:
 $$\begin{aligned}
 \mu^{3/5}\sup_{t\geq0} \|u(t)\|_{\dot B^{6/5}_{5/2,1}(\R^3)} + \|\mu\nabla^2u,u_t\|_{L_{5/2,1}(\R_+;L_{5/2}(\R^3))} &\leq C \mu^{3/5}\|u_0\|_{\dot B^{6/5}_{5/2,1}(\R^3)},\\
 \mu^{3/10}\sup_{t\geq0} \|u(t)\|_{\dot B^{3/5}_{10/7,1}(\R^3)} + \|\mu\nabla^2u,u_t\|_{L_{10/7,1}(\R_+;L_{10/7}(\R^3))} &\leq C \mu^{3/10}\|u_0\|_{\dot B^{3/5}_{10/7,1}(\R^3)},\\
  \mu^{7/10}\sup_{t\geq0} \|tu(t)\|_{\dot B^{7/5}_{10/3,1}(\R^3)} + \|\mu\nabla^2(tu),(tu)_t\|_{L_{10/3,1}(\R_+;L_{10/3}(\R^3))} &\leq C \mu^{7/10}\|u_0\|_{\dot B^{3/5}_{10/7,1}(\R^3)}\\
\andf\quad \mu \int_0^\infty\|\nabla u\|_{L_\infty(\R^3)}\,dx \leq C  \|u_0\|_{\dot B^{6/5}_{5/2,1}(\R^3)}^{1/3} &\|u_0\|_{\dot B^{3/5}_{10/7,1}(\R^3)}^{2/3}.
 \end{aligned}$$
\end{thm}
\begin{rmk} That $u_0$ belongs to $\dot B^{3/5}_{10/7,1}(\R^3)$ may be seen as a low frequency assumption
that exactly corresponds to the critical embedding in $L_2(\R^3)$ (hence our solutions have finite energies). Other choices of exponents might be possible both 
for low and high regularity. 

We also want to stress the fact that the smallness condition~\eqref{eq:smallness} is scaling invariant. 
\end{rmk}
The rest of the paper is structured as follows. 
In the next section, we prove our two-dimensional global result for~\eqref{eq:presless}. 
Section~\ref{s:NS3D} is devoted to the proof of Theorem~\ref{thm:NS3D}. 
In Appendix, we establish a maximal regularity estimate involving Lorentz spaces, of independent interest. 

We shall use standard notations and conventions. In particular, 
$C$ will always designate harmless constants that do not depend on `important' 
quantities, and we shall sometimes note $A\lesssim B$ instead of  $A\leq CB.$


\subsection*{Acknowledgments} This work was partially supported by   ANR-15-CE40-0011.

The second author (P.B.M.) has been also supported by the Polish National Science Centre's grant No 2018/30/M/ST1/00340  (HARMONIA).

We are extremely grateful  to the anonymous referee, for pointing out a number of severe mistakes in the first version of the paper.


 \section{The two dimensional case}\label{s:NS2D}

 This part is dedicated to the proof of Theorem~\ref{thm:pressureless}. 
 The key observation is that   combining the energy balance associated to~\eqref{eq:presless}
 with Proposition \ref{p:MR}  supplies a bound of the norm of 
 $$
\mbox{ $u$  in  $\dot W^{2,1}_{4/3,(4/3,1)}(\R^2\times\R^+)$ 
 and in $L_{4,1}(\R_+;L_4(\R^2))$ in terms of $u_0,$ }
 $$
 provided   the density is close to $1.$ 
 The latter information  will enable us to bound  $tu$ in $\dot W^{2,1}_{4,(4,1)}(\R^2\times\R^+)$
 and  to get a control on $\div u$ in $L_1(\R_+;L_\infty(\R^2)).$ Next, 
 reverting  to the mass conservation equation, one can ensure that  $\rho-1$ 
 remains small   provided   $\nu:=\mu+\mu'$ is large enough. 
  Then, we observe that the very same arguments leading to the control of $\div u$
 also allow to bound $\nabla u$ in $L_1(\R_+;L_\infty(\R^2)).$
   From that point,  we follow  the energy method of~\cite[Sec.~4]{DFP} 
  going to Lagrangian coordinates in order to prove uniqueness, 
  and the rigorous proof of existence is obtained by compactness arguments, after constructing a sequence 
  of smoother solutions (see the next section). 
  \medbreak
  Let us now go to the details of the proof.  In the first five steps, we assume  we are given a smooth enough
  solution $(\rho,u)$ of ~\eqref{eq:presless} and we concentrate on the proof of a priori estimates.
  We  suppose  that $\mu'\geq0$  and, to simplify the computations, we  take
 $\mu=1.$ That latter assumption is not restrictive, since $(\rho,u)$ satisfies~\eqref{eq:presless} 
 with coefficients $(\mu,\mu')$ if and only if 
 \begin{equation}\label{eq:rescaling}(\wt\rho,\wt u)(t,x):= (\rho, \mu^{-1}u)(\mu^{-1}t,x)\end{equation}
 satisfies~\eqref{eq:presless} with coefficients $1$ and $\mu'/\mu.$
  \smallbreak
  
 \subsubsection*{Step 1. The energy balance}
  As  already pointed out,  the space for $u_0$ is continuously embedded in   $L_2(\R^2).$
 Hence the initial data have finite energy.  Now, 
 the  energy balance for~\eqref{eq:presless}  reads
\begin{equation}
 \frac{1}{2} \frac{d}{dt} \int_{\R^2} \rho |u|^2 dx + 
 \int_{\R^2} \bigl( |\nabla u|^2 + \mu' (\div u)^2\bigr)dx =0.
\end{equation}
Provided the first part of~\eqref{eq:smallrho1} is fulfilled with small enough $c,$ we thus have, 
denoting  $\nu:=1+\mu',$  
\begin{equation}\label{eq:energy1} \|u\|_{L_\infty(\R_+;L_2(\R^2))}^2 + 
2\|\nabla\cP u\|_{L_2(\R_+;L_2(\R^2))}^2+2\nu\|\div u\|_{L_2(\R_+;L_2(\R^2))}^2
\leq  2\|u_0\|_{L_2}^2.
\end{equation}

 \subsubsection*{Step 2. Control of the norm of the solution in  $\dot W^{2,1}_{4/3,(4/3,1)}$}
 
Rewrite the velocity equation as:
 $$ u_t-\Delta u- \mu'\nabla\div u= f:=(1-\rho)u_t-\rho u\cdot\nabla u. $$
 Projecting the equation by means of  the Helmholtz projectors $\cP$ and $\cQ,$  we get
 $$ (\cP u)_t-\Delta\cP u = \cP f\andf (\cQ u)_t-\nu\Delta\cQ u=\cQ f.$$
 Applying  Proposition \ref{p:MR} with $d=2,$ $p=q=4/3,$  $m=s=4$ and $r=1$ yields 
 \begin{multline}\label{sss1}
\|\cP u\|_{L_\infty(\R_+;\dot B^{1/2}_{4/3,1}(\R^2))} 
+  \|\nabla^2\cP u,\cP u_t\|_{L_{4/3,1}(\R_+\!;L_{4/3}(\R^2))} +\|\cP u\|_{L_{4,1}(\R_+;L_4(\R^2))}
\\\lesssim
\|\cP u_0\|_{\dot B^{1/2}_{4/3,1}(\R^2)} + 
\|\cP f\|_{L_{4/3,1}(\R_+\!;L_{4/3}(\R^2))}\end{multline}
and\footnote{Observe that we can also write $\nabla^2\cQ u$ instead of $\nabla \div u$  in~\eqref{sss2} and that $\nabla \div\cQ u = \nabla \div u$.}
\begin{multline}\label{sss2}
 \nu^{1/4}\|\cQ u\|_{L_\infty(\R_+;\dot B^{1/2}_{4/3,1} (\R^2))}
+  \|\nu\nabla\div u, \cQ u_t\|_{L_{4/3,1}(\R_+;L_{4/3}(\R^2))} 
 +\nu^{1/2} \|\cQ u\|_{L_{4,1}(\R_+;L_4(\R^2))} \\
\lesssim \nu^{1/4}\|\cQ u_0\|_{\dot B^{1/2}_{4/3,1}(\R^2)} + 
\|\cQ f\|_{L_{4/3,1}(\R_+;L_{4/3}(\R^2))}.
\end{multline}
To bound  $\cP f$ and $\cQ f,$ we use the fact that 
  $\cP$ and $\cQ$ are continuous on  $L_{4/3,1}(\R_+ ; L_{4/3} (\R^2)),$ so that it  is  enough  to 
estimate $f$ in $L_{4/3,1}(\R_+;L_{4/3}(\R^2)).$  Now, using H\"older inequality in Lorentz spaces (see \eqref{eq:holder}), we find that
$$\displaylines{
\|f\|_{L_{4/3,1}(\R_+;L_{4/3}(\R^2))} \leq \|1-\rho\|_{L_\infty(\R_+\times\R^2)} \|u_t\|_{L_{4/3,1}(\R_+;L_{4/3}(\R^2))} 
\hfill\cr\hfill+\|\rho\|_{L_\infty(\R_+\times\R^2)}\|\nabla u\|_{L_{2}(\R_+\times\R^2)}
\|u\|_{L_{4,1}(\R_+;L_{4}(\R^2))}.}$$
Hence, summing up \eqref{sss1} and \eqref{sss2}, and assuming smallness of $\rho -1$,  one gets
$$\displaylines{\|\cP u\|_{L_\infty(\R_+;\dot B^{1/2}_{4/3,1} (\R^2))} + \nu^{1/4}\|\cQ u\|_{L_\infty(\R_+;\dot B^{1/2}_{4/3,1}(\R^2))}+
\|u_t\|_{L_{4/3,1}(\R_+;L_{4/3}(\R^2))} \hfill\cr\hfill+  \|\nabla^2\cP u,\nu\nabla\div u\|_{L_{4/3,1}(\R_+;L_{4/3}(\R^2))} 
 +\|\cP u\|_{L_{4,1}(\R_+;L_4(\R^2))}+ \nu^{1/2} \|\cQ u\|_{L_{4,1}(\R_+;L_4(\R^2))}
 \hfill\cr\hfill\lesssim \|\cP u_0\|_{\dot B^{1/2}_{4/3,1}(\R^2)} + \nu^{1/4}\|\cQ u_0\|_{\dot B^{1/2}_{4/3,1}(\R^2)} + 
 \|\nabla u\|_{L_{2}(\R_+\times\R^2)}
\|u\|_{L_{4,1}(\R_+;L_{4}(\R^2))}.}$$
 In the case $\nu\geq1,$ it is now clear that there exists a constant $\eta>0$ independent of $\nu$ such that, if 
 \begin{equation}\label{eq:smallDu}
   \|\nabla u\|_{L_{2}(\R_+\times\R^2)}\leq \eta,
   \end{equation}
   then we eventually have 
$$\displaylines{\|\cP u\|_{L_\infty(\R_+;\dot B^{1/2}_{4/3,1}(\R^2))} + \nu^{1/4}\|\cQ u\|_{L_\infty(\R_+;\dot B^{1/2}_{4/3,1}(\R^2))}+
\|u_t\|_{L_{4/3,1}(\R_+;L_{4/3}(\R^2))} \hfill\cr\hfill+  \|\nabla^2\cP u,\nu\nabla\div u\|_{L_{4/3,1}(\R_+;L_{4/3}(\R^2))} 
 +\|\cP u\|_{L_{4,1}(\R_+;L_4(\R^2))}+ \nu^{1/2} \|\cQ u\|_{L_{4,1}(\R_+;L_4(\R^2))}
 \hfill\cr\hfill\lesssim \|\cP u_0\|_{\dot B^{1/2}_{4/3,1}(\R^2)} + \nu^{1/4}\|\cQ u_0\|_{\dot B^{1/2}_{4/3,1}(\R^2)}.}$$
If \eqref{eq:smallDu} is not satisfied, then  the idea is to split the time line into a finite number $K$ of 
intervals $[T_{k-1},T_{k})$ with $T_0 = 0$ and $T_K=\infty,$ 
 $$\| \nabla u \|_{L_2 ((T_{k - 1} , T_k) \times \R^2)} = \eta\ \hbox{ for all }\ 1 \leq k \leq K - 1,\andf\| \nabla u \|_{L_2 ((T_{K - 1} , T_K) \times \R^2)} \leq \eta. $$
For given $\eta$, we now calculate how many time intervals are needed for this splitting. We have 
$$ K\eta^2\geq  \sum_{k = 1}^{K } \| \nabla u \|_{L_2 ((T_{k - 1} , T_k) \times \R^2)}^2 = \| \nabla u \|_{L_2 (\R_+ \times \R^2)}^2 \geq \sum_{k = 1}^{K - 1} \| \nabla u \|_{L_2 ((T_{k - 1} , T_k) \times \R^2)}^2 = (K - 1) \eta^2,$$
whence  
\begin{equation}\label{sss0}
K=\lceil \eta^{-2} \| \nabla u \|_{L_2 (\R_+ \times \R^2)}^2\rceil.
\end{equation}
Having this information at hand, we adapt \eqref{sss1} and \eqref{sss2} to the finite time interval setting, getting, for each $k\in\{0,\cdots K-1\},$ 
\begin{multline}\label{eq:step2x}\|\cP u\|_{L_\infty(T_k,T_{k+1};\dot B^{1/2}_{4/3,1}(\R^2))}
    + \nu^{1/4}\|\cQ u\|_{L_\infty(T_k,T_{k+1};\dot B^{1/2}_{4/3,1}(\R^2))} 
   +\|\cP u\|_{L_{4,1}(T_k,T_{k+1};L_4(\R^2))}\\+ \nu^{1/2} \|\cQ u\|_{L_{4,1}(T_k,T_{k+1};L_4(\R^2))}
    + \|u_t\|_{L_{4/3,1}(T_k,T_{k+1};L_{4/3}(\R^2))}
+  \|\nabla^2 \cP u\|_{L_{4/3,1}(T_k,T_{k+1};L_{4/3}(\R^2))} 
\\+ \nu \|\nabla\div u\|_{L_{4/3,1}(T_k,T_{k+1};L_{4/3}(\R^2))}
\leq C\big(\|\cP u(T_k)\|_{\dot B^{1/2}_{4/3,1}(\R^2)} + \nu^{1/4}\|\cQ u(T_k)\|_{\dot B^{1/2}_{4/3,1}(\R^2)}\big)\cdotp
\end{multline}
Arguing by induction  and remembering that  $K$ is estimated by  \eqref{sss0} and \eqref{eq:energy1}, we conclude that 
\begin{multline}\label{eq:L41}
\|\cP u\|_{L_\infty(\R_+;\dot B^{1/2}_{4/3,1}(\R^2))}+\nu^{1/4}\|\cQ u\|_{L_\infty(\R_+;\dot B^{1/2}_{4/3,1}(\R^2))} 
 + \|\cP u\|_{L_{4,1}(\R_+;L_4(\R^2))} \\ +\nu^{1/2}\|\cQ  u\|_{L_{4,1}(\R_+;L_4(\R^2))}
    + \|u_t\|_{L_{4/3,1}(\R_+;L_{4/3}(\R^2))} 
+  \|\nabla^2\cP u\|_{L_{4/3,1}( \R_+;L_{4/3}(\R^2))} \\
+ \nu \|\nabla\div u\|_{L_{4/3,1}( \R_+;L_{4/3}(\R^2))} 
\leq 
 C\big( \|\cP u_0\|_{\dot B^{1/2}_{4/3,1}(\R^2)} + 
\nu^{1/4}\|\cQ u_0\|_{\dot B^{1/2}_{4/3,1}(\R^2)} \big)e^{C\|u_0\|_{L_2(\R^2)}^2}.
\end{multline}

  \subsubsection*{Step 3. A time weighted estimate}

We now look at the momentum equation in the form
\begin{equation}\label{t-NS}
  \begin{array}{l}
  (t\, u)_t   -  \Delta (tu) - \mu'\nabla\div (tu)=(1-\rho)(t\,u)_t+\rho u- t \rho u \cdot \nabla u.
  \end{array}
 \end{equation}

 By definition, the initial datum for $tu$ is  zero, and~\eqref{eq:L41} provides us with a bound  
for  the term $\rho u$  in  the space $L_{4,1}(\R_+;L_4(\R^2)).$ This  
gives us some hint on the regularity of the whole right-hand side. 
Now, projecting~\eqref{t-NS} by means of  $\cP$ and $\cQ,$ 
using the  maximal regularity estimates of 
Proposition \ref{p:MR}, and still assuming the first part of~\eqref{eq:smallrho1}, 
one gets for all $0\leq T\leq T'\leq\infty,$ 
 \begin{multline}\label{eq:g0}
 \sup_{T\leq t\leq T'} \|t\cP u\|_{\dot B^{3/2}_{4,1}(\R^2)} 
 + \nu^{3/4}\sup_{T\leq t\leq T'} \|t\cQ u\|_{\dot B^{3/2}_{4,1}(\R^2)} + 
\|\nabla^2 t\cP u\|_{L_{4,1}(T,T';L_4(\R^2))} \\
+ \|(tu)_t\|_{L_{4,1}(T,T';L_4(\R^2))} 
+ \nu \|\nabla\div tu\|_{L_{4,1}(T,T';L_4(\R^2))} 
\leq C\Bigl(  \| T \cP u(T)\|_{\dot B^{3/2}_{4,1}(\R^2)}\\ + \nu^{3/4} \| T\cQ u(T)\|_{\dot B^{3/2}_{4,1}(\R^2)} + 
 \| u \|_{L_{4,1}(T,T';L_4(\R^2))} 
+ \|tu\cdot\nabla u\|_{L_{4,1}(T,T';L_4(\R^2))}\Bigr)\cdotp\end{multline}
Since we have $\dot B^{1/2}_{4,1}(\R^2)\hookrightarrow L_\infty(\R^2),$ the term with $u\cdot\nabla u$ may be bounded as follows:
\begin{equation}\label{w8}
 \|tu\cdot\nabla u\|_{L_{4,1}(T,T';L_4(\R^2))} 
 \leq C\|u\|_{L_{4,1}(T,T';L_4(\R^2))} \|t  u\|_{L_\infty(T,T'; \dot B^{3/2}_{4,1}(\R^2))}.\end{equation}
Consequently, if $\nu \geq 1$ there exists a constant $\eta>0$ independent of $\nu$ such that, if 
\begin{equation}\label{eq:key}
 \|u\|_{L_{4,1}(T , T' ; L_4 (\R^2))}\leq \eta,
 \end{equation}
then Inequality~\eqref{eq:g0} reduces to 
 \begin{multline}\label{eq:g1}
  \sup_{T\leq t\leq T'} \|t\cP u\|_{\dot B^{3/2}_{4,1}(\R^2)} 
 + \nu^{3/4}\sup_{T\leq t\leq T'} \|t\cQ u\|_{\dot B^{3/2}_{4,1}(\R^2)} + 
\|\nabla^2 t\cP u\|_{L_{4,1}(T,T';L_4(\R^2))} \\
+ \|(tu)_t\|_{L_{4,1}(T,T';L_4(\R^2))} 
+ \nu \|\nabla\div tu\|_{L_{4,1}(T,T';L_4(\R^2))} 
\leq C\Bigl(  \|T \cP u(T)\|_{\dot B^{3/2}_{4,1}(\R^2)}\\ + \nu^{3/4} \| T\cQ u(T)\|_{\dot B^{3/2}_{4,1}(\R^2)} + 
 \| u \|_{L_{4,1}(T,T';L_4(\R^2))} \Bigr)\cdotp\end{multline}
Clearly, if one can take $T=0$ and $T'=\infty$ in~\eqref{eq:key}, then we 
control the left-hand side of~\eqref{eq:g0} on $\R_+,$  so 
assume from now on that 
$ \|u\|_{L_{4,1}(\R_+;L_4(\R^2))}>\eta.$
We claim that there exists  a finite sequence
$0=T_0<T_1<\cdots< T_{K-1}< T_{K}=\infty$
such that \eqref{eq:key} if fulfilled on $[T_k,T_{k+1}]$ for each $k\in\{0,\cdots, K-1\}.$
In order to prove our claim, we  introduce 
$$ U(t):=\|u(t,\cdot)\|_{L_{4}(\R^2)}$$ and recall that
\begin{equation}
 \|U\|_{L_{4,1}(\R_+)} = 4 \int_0^\infty |\{t\in \R_+: |U(t)|>s\}|^{1/4}\, ds.
\end{equation}
From Lebesgue dominated convergence theorem,   we have
\begin{equation}\label{w7b} \|U\|_{L_{4,1}(T',T'')} = 4 \int_0^\infty |\{t\in (T',T''): |U(t)|>s\}|^{1/4} ds \to 0\quad\hbox{as}\quad
 T''-T'\to0.\end{equation}
 Hence one can construct inductively a family  $0=T_0<T_1<\dotsm<T_{k}<\dotsm < T_K = \infty$
 such that
\begin{equation}\label{eq:U}
 \|U\|_{L_{4,1}(T_{k - 1},T_{k})}=\eta \quad \text{for} \quad 1 \leq k \leq K - 1 \andf\|U\|_{L_{4,1}(T_{K - 1},T_{K})} \leq \eta.
\end{equation}
By simple H\"older inequality on series, we have
$$
 \sum_{k=1}^{K - 1} | \{ t\in (T_{k-1},T_k): |U(t)|>s\} |^{1/4} \leq 
 (K - 1)^{3/4}  \biggl( \sum_{k=1}^{K - 1} |\{ t \in (T_{k-1},T_k): |U(t)|>s\}|\biggr)^{1/4}.$$
Hence, integrating with respect to $s,$ and using \eqref{eq:U}, we get
$$(K-1)\eta\leq (K-1)^{3/4}\int_0^\infty 
|\{t\in(0,T_{K-1}): |U(t)|>s\}|^{1/4}\,ds,
$$
whence 
$$(K-1)^{1/4}\eta\leq \|U\|_{L_{4,1}(\R_+)}.$$
Now, remembering the definition of $U$ and \eqref{eq:L41}, one may conclude that
there exists some constant $C$ such that 
\begin{equation}\label{eq:K}
K=[ C  \bigl(\|\cP u_0\|_{\dot B^{1/2}_{4/3,1}(\R^2)}^4 + 
\nu\|\cQ u_0\|_{\dot B^{1/2}_{4/3,1}(\R^2)}^4\bigr) e^{C\|u_0\|_{L_2(\R^2)}^2}]+1.
\end{equation}
Let $X_k:=\sup_{t\in (T_k,T_{k+1})} \bigl(\|t\cP u(t)\|_{\dot B^{3/2}_{4,1}(\R^2)}
+ \nu^{3/4} \|t\cQ u(t)\|_{\dot B^{3/2}_{4,1}(\R^2)}\bigr)$.
Then,  \eqref{eq:g1} and \eqref{eq:U} ensure that  $X_0\leq C\eta$ and that
\begin{equation}
 X_k \leq C(\eta + X_{k-1})\quad\hbox{for all }\ k\in\{1,\cdots,K-1\}.
\end{equation}
So,  arguing by induction, we  eventually get for all $m\in\{0,\cdots,K-1\},$
\begin{equation}\label{eq:g2}
 X_m \leq C\eta \sum_{\ell=0}^mC^\ell  \leq  \frac{C^{K+1}}{C-1}\:\eta.
\end{equation}
Reverting to \eqref{eq:g1}, then using  \eqref{eq:K}, we conclude that
for some universal constant $C\geq1,$ 
\begin{multline}\label{w12}
 \sup_{t\geq0} \|t\cP u\|_{\dot B^{3/2}_{4,1}(\R^2)} + 
 \nu^{1/4} \sup_{t\geq0} \|t\cQ u\|_{\dot B^{3/2}_{4,1}(\R^2)} +
 \|\nabla^2  t\cP u\|_{L_{4,1}(\R_+;L_4(\R^2))}\\ 
+ \|(tu)_t\|_{L_{4,1}(\R_+;L_4(\R^2))} 
+ \nu \|\nabla\div tu\|_{L_{4,1}(\R_+;L_4(\R^2))}\\ 
\leq C \exp\biggl(C \bigl(\|\cP u_0\|_{\dot B^{1/2}_{4/3,1}(\R^2)}^4 + 
\nu\|\cQ u_0\|_{\dot B^{1/2}_{4/3,1}(\R^2)}^4\bigr) e^{C\|u_0\|_{L_2(\R^2)}^2}\biggr)\cdotp \end{multline}

  \subsubsection*{Step 4. Bounding $\div u$}

In order to keep the density close to $1,$ we have to bound 
 $\div u$  in $L_1(\R_+;L_\infty(\R^2)).$  The key observation is that
\begin{equation}\label{c2}
 \div (t u) \in L_{4,1}(\R_+;\dot W^1_{4}(\R^2)) 
 \andf\div u \in  L_{4/3,1}(\R_+;\dot W^{1}_{4/3}(\R^2)).
\end{equation}
Now, from Gagliardo-Nirenberg inequality and Sobolev embedding, we see that 
\begin{equation}\label{c3}
 \|z\|_{L_\infty(\R^2)} \leq C \|\nabla z\|^{1/2}_{L_4(\R^2)} \|\nabla z\|^{1/2}_{L_{4/3}(\R^2)}.
\end{equation}
So we have, thanks to H\"older inequality in 
Lorentz spaces\footnote{Here it is essential to have a  Lorentz space \emph{with last index $1,$ as regards  time exponent.}}, 
$$\begin{aligned}
\int_0^\infty \|\div u\|_{L_\infty}\, dt  &
\leq C \int_0^\infty t^{-1/2} \|t\nabla\div u\|_{L_4}^{1/2} 
\|\nabla \div u\|^{1/2}_{L_{4/3}}\,dt  \\ &\leq C \|t^{-1/2}\|_{L_{2,\infty}(\R_+)}
\bigl\|\|t\nabla\div u\|_{L_{4}(\R^2)}^{1/2}\bigr\|_{L_{8,2}(\R_+)}
\bigl\|\|\nabla\div u\|_{L_{4/3}(\R^2)}^{1/2}\bigr\|_{L_{8/3,2}(\R_+)}\\
 &\leq C \|t\nabla\div u\|_{L_{4,1}(\R_+;L_4(\R^2))}^{1/2}
\|\nabla\div u\|_{L_{4/3,1}(\R_+;L_{4/3}(\R^2))}^{1/2}. \end{aligned}
$$
Hence, thanks to \eqref{eq:L41} and \eqref{w12}, 
\begin{multline}\label{eq:step4}
 \nu\int_0^\infty \|\div u\|_{L_\infty}\, dt \leq C\bigl(\|\cP u_0\|_{\dot B^{1/2}_{4/3,1}(\R^2)}^{1/2} + 
\nu^{1/8}\|\cQ u_0\|_{\dot B^{1/2}_{4/3,1}(\R^2)}^{1/2}\bigr) \\\times
 \exp\biggl(C \bigl(\|\cP u_0\|_{\dot B^{1/2}_{4/3,1}(\R^2)}^4 + 
\nu\|\cQ u_0\|_{\dot B^{1/2}_{4/3,1}(\R^2)}^4\bigr) e^{C\|u_0\|_{L_2(\R^2)}^2}\biggr)\cdotp
 \end{multline}

  \subsubsection*{Step 5. Bounding the density} 

The discrepancy of the density to $1$ (that is $a:=\rho-1$) may be controlled by means of the mass equation: 
$$\partial_ta +u\cdot\nabla a + (1+a)\div u=0,$$
which gives 
$$\|a(t)\|_{L_\infty(\R^2)}\leq \|a_0\|_{L_\infty(\R^2)} + \int_0^t\|\div u\|_{L_\infty(\R^2)}\,d\tau 
+ \int_0^t\|a\|_{L_\infty(\R^2)}\|\div u\|_{L_\infty(\R^2)}\,d\tau,$$
and thus, owing to the integral version of Gronwall lemma,
$$\begin{aligned}
\|a(t)\|_{L_{\infty}(\R^2)} &\leq \bigg( \| a_0 \|_{L_{\infty}(\R^2)} + \int_0^t
 e^{-\int_0^s \| {\rm div}\, u \|_{L_{\infty}(\R^2)}\,ds}
 \|  \div u\|_{L_{\infty}(\R^2)}\, d\tau \bigg) e^{\int_0^t \| {\rm div}\, u \|_{L_{\infty}(\R^2)}\,dt}\\[5pt]
 &= \displaystyle 
 \| a_0 \|_{L_{\infty}(\R^2)}e^{\int_0^t \| {\rm div}\, u \|_{L_{\infty}(\R^2)}\,dt}\;
 + \;e^{\int_0^t \| {\rm div}\, u \|_{L_{\infty}(\R^2)}\,dt}-1.\end{aligned}$$
 Hence, provided we have \eqref{eq:smallrho0} and $c \leq 1$ and
\begin{equation}\label{eq:condsmall}
\int_0^T\|\div u\|_{L_\infty(\R^2)}\,d\tau\leq\log(1+c/2),
\end{equation}
the smallness property~\eqref{eq:smallrho1} is satisfied on $[0,T].$
Bearing in mind Inequality~\eqref{eq:step4} and the  assumption~\eqref{eq:data}, we get
$$ \nu \int_0^{\infty} \| \div u \|_{L_{\infty}(\R^2)} \, d t \leq C  M^{\frac{1}{2}} e^{C M^4}.$$
Consequently, in order to have  \eqref{eq:condsmall}  satisfied, it suffices that
$$
\nu^{-1} C  M^{\frac{1}{2}} e^{C M^4}\leq \log(1+c/2),$$
which corresponds to the condition on $\nu$ given in  the statement of the theorem. 

\subsubsection*{Step 6: Uniqueness} 

The key to uniqueness is  that $\nabla u$ is in $L_1(\R_+;L_\infty(\R^2)).$
To get that property, one may proceed
exactly as for bounding $\div u,$ writing that
\begin{equation}\label{u-1}
\int_0^\infty \|\nabla u\|_{L_\infty(\R^2)}\, dt \leq C \|t^{-1/2}\|_{L_{2,\infty}(\R_+)}
\|t\nabla^2u\|_{L_{4,1}(\R_+;L_4(\R^2))}^{1/2}\|\nabla^2 u\|_{L_{4/3,1}(\R_+;L_{4/3}(\R^2))}^{1/2},
\end{equation}
and using that the right-hand side is bounded in terms of $u_0$ 
according to~\eqref{eq:L41} and~\eqref{w12}. 
\medbreak
Because of the hyperbolic nature of the continuity equation, the uniqueness issue is not straightforward, as the regularity of the density is very low. 
However, knowing that  $\nabla u$ is in $L_1(\R_+;L_\infty(\R^2))$ 
enables us  to rewrite our system in Lagrangian coordinates. 
More precisely, for all $y\in\R^2,$ consider the following ODE:
\begin{equation}\label{eq:ODE}
 \frac{dX}{dt}(t,y)=u(t,X(t,y)), \qquad \qquad X|_{t=0}=y.
\end{equation}
Having~\eqref{u-1} at hand guarantees that~\eqref{eq:ODE} 
defines a $C^1$ flow $X$ on $\R_+\times\R^2.$ 
\medbreak
Let us express the density and velocity in the  new coordinates:
\begin{equation}
 \eta(t,y)=\rho(t,X(t,y)), \qquad v(t,y)=u(t,X(t,y)).
\end{equation}
Then, the  system for $(\eta,v)$ reads (see  details  in, e.g.,~\cite{D-Fourier}): 
\begin{equation}\label{u1}
 \begin{array}{l}
  (J_v\eta)_t=0, \\[5pt]
  \rho_0 v_t -\divv(\nabla_vv +  \mu' (\divv v)\Id) =0,
 \end{array}
\end{equation}
where $\nabla_v:=A_v^\top \nabla_y$, 
$\div_v:= \div(J_v^{-1}A_v\cdot)= A_v^\top:\nabla_y$ with $A_v=(DX_v)^{-1}$ and $J_v=\det(DX_v).$
One points out that $J_v^{-1}A_v=\adj(DX_v)$ (the adjugate matrix of $DX_v$). 
\medbreak
Since, in our framework  
the Lagrangian and Eulerian formulations are equivalent (see, e.g.,~\cite{D-Fourier,Mu03}), 
it suffices to prove uniqueness at the level of Lagrangian coordinates. 
Therefore, consider two  solutions $(\eta,v)$ and $(\bar \eta, \bar v)$ of
\eqref{u1} emanating from the data $(\rho_0,u_0).$ Then, the  difference 
of velocities  $\dv: = \bv -  v$
satisfies  
\begin{multline}\label{eq:dv}
 \rho_0 \dv_t -\divv(\nabla_v\dv +  \mu' (\divv\dv)\Id) \\
=\bigl(\div_{\!\bv}\nabla_\bv-\div_v\nabla_v\bigr)\bv +\mu'\bigl(\div_{\!\bv}\Id\div_{\!\bv}-\divv\Id\divv\bigr)\bv.
\end{multline}
Note that 
$$\begin{aligned}
\bigl(\div_{\!\bv}\nabla_\bv-\div_v\nabla_v\bigr)\bv&=\div\bigl((\adj(DX_\bv)A_\bv^\top-
\adj(DX_v)A_v^\top)\cdot\nabla\bv\bigr),\\
\bigl(\div{\!\bv}\Id\div_{\!\bv}-\divv \Id\divv\bigr)\bv&= \div\bigl((\adj(DX_\bv)A_\bv^\top-
\adj(DX_v)A_v^\top):\nabla\bv\bigr)\cdotp
\end{aligned}$$
Now,  taking the $L^2$ scalar product of~\eqref{eq:dv} with $\dv$ and integrating by parts delivers
$$\displaylines{\quad\frac12\frac d{dt}\|\sqrt{\rho_0}\,\dv\|_{L_2(\R^2)}^2
+\|\nabla_v\dv\|_{L_2(\R^2)}^2+\mu'\|\divv\dv\|_{L_2(\R^2)}^2\hfill\cr\hfill\leq
\nu\bigl\|\bigl(\adj(DX_\bv)A_\bv^\top-
\adj(DX_v)A_v^\top\bigr)\cdot\nabla\bv\bigr\|_{L^2(\R^2)}\|\nabla\dv\|_{L_2(\R^2)}.\quad}
$$
Let us take an interval $[0,T]$ for which 
\begin{equation}\label{eq:issmall}
\max\biggl(\int_0^T\|\nabla v\|_{L_\infty(\R^2)}\,dt,\int_0^T\|\nabla \bv\|_{L_\infty(\R^2)}\,dt\biggr)
\quad\hbox{is small.}
\end{equation}
All terms like $\Id-A_w,$ $1-J_w$  or $\Id-\adj(DX_w)$ (with $w=v,\bv$) 
may be computed by using Neumann series expansions, and  
we end up with  pointwise estimates of the following type:
$$ |\Id - A_w|\lesssim \int_0^t |\nabla w|\, dt', \qquad 
 |\Id-\adj(DX_w)|\lesssim \int_0^t |\nabla  w| \, dt', \qquad
 |1-J_w|\lesssim \int_0^t|\nabla w| \, dt'.$$
{}From this, we deduce that 
$$\frac d{dt}\|\sqrt{\rho_0}\,\dv\|_{L_2(\R^2)}^2
+\|\nabla\dv\|_{L_2(\R^2)}^2\lesssim (\|\nabla v\|_{L_\infty(\R^2)}\!+\!\|\nabla\bv\|_{L_\infty(\R^2)})
\|\nabla\dv\|_{L_2(\R^2)}\biggl\|\int_0^t \nabla\dv\,d\tau\biggr\|_{L_2(\R^2)}. 
$$
Because we have, by Cauchy-Schwarz inequality, 
$$t^{-1/2} \biggl\|\int_0^t \nabla\dv\,d\tau\biggr\|_{L_2(\R^2)} \leq \|\nabla\dv\|_{L_2((0,t)\times\R^2)},$$
integrating the above inequality (and using again Cauchy-Schwarz inequality) yields
$$\|\sqrt{\rho_0}\,\dv(t)\|_{L_2(\R^2)}^2+\int_0^t \|\nabla\dv\|_{L_2(\R^2)}^2\,d\tau\leq
C\biggl(\int_0^t \tau\|(\nabla v,\nabla\bar v)(\tau)\|_{L_\infty(\R^2)}^2d\tau\biggr)^{1/2}
\int_0^t \|\nabla\dv\|_{L_2(\R^2)}^2\,d\tau.$$
Hence,  there exists $c>0$ such that if, in addition to~\eqref{eq:issmall}, 
we have
\begin{equation}\label{eq:issmall2}
\|t^{1/2}\nabla w\|_{L_2(0,T;L_\infty(\R^2))}\leq c\quad\hbox{for }\ w=v,\bv, 
\end{equation}
then we have $\dv\equiv0$ on $[0,T],$ that is to say $\bv=v.$
Since $\deta=(J_{\bv}^{-1}-J_v)\rho_0,$  we get  $\bar\eta=\eta,$ too. 
\medbreak
In light of the above arguments, in order to get uniqueness on the whole $\R_+,$
it suffices to show that our solutions satisfy not only $\nabla u\in L_1(\R_+;L_\infty(\R^2)),$ 
but also 
\begin{equation}\label{eq:safdasdfafds}
\int_0^\infty t\|\nabla u\|_{L_\infty(\R^2)}^2\,dt<\infty.
\end{equation} 
This is a consequence of~\eqref{c3}, as it gives
$$\begin{aligned}
\int_0^\infty t\|\nabla u\|_{L_\infty(\R^2)}^2\,dt&\lesssim
\int_0^\infty \|t\nabla^2u\|_{L_4(\R^2)}\|\nabla^2u\|_{L_{4/3}(\R^2)} \, d t\\
&\lesssim \|t\nabla^2u\|_{L_4(\R_+\times\R^2)} \|\nabla^2 u\|_{L_{4/3}(\R_+\times\R^2)}\\
&\lesssim \|t\nabla^2u\|_{L_{4,1}(\R_+;L_4(\R^2))} \|\nabla^2 u\|_{L_{4/3,1}(\R_+;L_{4/3}(\R^2))}.
\end{aligned}
$$

\subsubsection*{Step 7: Proof of existence}
The idea is to smooth out the data,  and to 
solve~\eqref{eq:presless} supplemented with those data, 
according to the local-in-time existence result of~\cite{DFP}
(that just requires the initial velocity to be smooth enough, and the initial density
to be close to $1$ in $L_\infty$). 
Then, the previous steps provide uniform bounds
that allow to show that those smoother solutions are actually global, 
 and  one can eventually pass to the limit. 
The reader may refer to the end of the next part
where more details  are given  on the existence issue, both for Theorems~\ref{thm:pressureless}
and~\ref{thm:NS3D}.

\section{The three dimensional case}\label{s:NS3D}

Our aim here is to prove  a global existence result  in the small data regime case
for System~\eqref{eq:NS3D}. In order to get the optimal dependency of the smallness condition in terms of the viscosity 
coefficient,  it is wise to resort  again  to the rescaling~\eqref{eq:rescaling}. 
So we assume from now on that  $\mu=1.$

The bulk of the proof consists in exhibiting global-in-time bounds in terms of the data for 
\begin{eqnarray}\label{eq:Xi}\Xi&\!\!\!:=\!\!\!&  \sup_{t\geq0} \| u(t) \|_{\dot B^{6/5}_{5/2,1}(\R^3)} + 
 \|\nabla^2u, u_t\|_{L_{5/2,1}(\R_+;L_{5/2}(\R^3))}\\\label{eq:Psi}\andf
 \Psi&\!\!\!:=\!\!\!&  \sup_{t\geq0} \| u(t) \|_{\dot B^{3/5}_{10/7,1}(\R^3)} + 
 \|\nabla^2u, u_t\|_{L_{10/7,1}(\R_+; L_{10/7}(\R^3))}.\end{eqnarray}
 {}From that control and Proposition \ref{p:MR}, we will estimate $u$ 
 in the   space $L_{10/3,1}(\R_+; L_{10/3}(\R^3))$
 (that will play  the same role as $L_{4,1}(\R_+; L_4(\R^2))$ 
 for \eqref{eq:presless}), 
 and exhibit a bound for $tu$ in the space $\dot W^{2,1}_{10/3,(10/3,1)}(\R_+ \times \R^3).$
 This will  eventually enable us to bound 
 $\nabla u$ in $L_1(\R_+;L_\infty(\R^3)).$
 {}From that stage, the proof of uniqueness follows the lines of  the two-dimensional case. 
 \medbreak
 \subsubsection*{Step 1. Control by the energy} Remembering that 
our assumptions  imply that $u_0$ is in $L_2(\R^3),$ 
we start with  the basic  energy balance:
\begin{equation}\label{eq:energy}
 \frac{d}{dt} \int_{\R^3} \rho |u|^2 \,dx + \int_{\R^3} |\nabla u|^2\, dx =0.
\end{equation}
By Sobolev embedding and  provided that \eqref{eq:smallrho} is satisfied,
this implies the following bound on $u$:
\begin{equation}\label{eq:energy3}
\|u\|_{L_\infty(\R_+;L_2(\R^3))} + \|u\|_{L_2(\R_+;\dot H^1(\R^3))}\lesssim \|u_0\|_{L_2(\R^3)}.  
\end{equation}
That relation will  enable us to control higher norms of the
solution, globally in time,  provided  some scaling invariant quantity
involving $u_0$ is small  enough. 

 \subsubsection*{Step 2. Control of the high norm} 
This  step is somehow standard: we want to construct smooth solutions like for the 
classical Navier-Stokes system. 
Now, assuming~\eqref{eq:smallrho} and  taking advantage of 
the maximal regularity estimate for the heat equation in 
$L_{5/2,1}(\R_+ ; L_{5/2}(\R^3))$ stated in Proposition~\ref{p:MR} yields (recall the definition of $\Xi$ in~\eqref{eq:Xi}):
$$\Xi \leq C\bigl(\|u\cdot\nabla u\|_{L_{5/2,1}(\R_+ ; L_{5/2}(\R^3))} 
 +\|u_0\|_{\dot B^{6/5}_{5/2,1}(\R^3)}\bigr)\cdotp$$
We see by H\"older inequality and Sobolev embedding $\dot W^1_{5/2}(\R^3) 
\hookrightarrow L_{15}(\R^3)$ that 
$$\begin{aligned}
 \|u\cdot \nabla u \|_{L_{5/2,1}(\R_+; L_{5/2}(\R^3))} &\leq C\|u\|_{L_\infty(\R_+;L_{3}(\R^3))} 
 \|\nabla u \|_{L_{5/2,1}(\R_+;L_{15}(\R^3))}\\ &\leq C\|u\|_{L_\infty(\R_+;L_{3}(\R^3))}
 \|\nabla^2u\|_{L_{5/2,1}(\R_+;L_{5/2}(\R^3))}.
\end{aligned}$$
Moreover,  we note that  by H\"older inequality, Sobolev embedding  and~\eqref{eq:energy3}, we have
\begin{equation}\label{eq:L31}
 \|u\|_{L_{3}(\R^3)} \leq C\|u\|_{L_2(\R^3)}^{2/3}\|u\|_{L_\infty(\R^3)}^{1/3}\leq 
 C\|u_0\|^{2/3}_{L_2(\R^3)} \|u\|^{1/3}
 _{\dot B^{6/5}_{5/2,1}(\R^3)}.
\end{equation}
Hence, altogether, this gives 
$$ \Xi \leq C\bigl(\|u_0\|_{L_2(\R^3)}^{2/3}\Xi^{1+1/3} +\Xi_0\bigr)\cdotp
$$
From this, we deduce that
\begin{equation}\label{t9}
(2C)^{4/3}\Xi_0^{1/3} \|u_0\|_{L_2}^{2/3} \leq 1
\quad\hbox{implies}\quad \Xi\leq 2 C \Xi_0. 
\end{equation}

 \subsubsection*{Step 3. Control of the low norm}  It is now a matter of bounding
 the functional $\Psi$ defined in~\eqref{eq:Psi}.  Thanks to Proposition \ref{p:MR}, we have
 \begin{equation}\label{t10}
   \Psi\leq C\bigl(\|u\cdot\nabla u\|_{L_{10/7,1}(\R_+; L_{10/7}(\R^3))} 
 +\|u_0\|_{\dot B^{3/5}_{10/7,1}(\R^3)}\bigr)\cdotp\end{equation}
 By H\"older inequality and Sobolev embedding $\dot W^1_{10/7}(\R^3) 
\hookrightarrow L_{30/11}(\R^3),$  we discover that
$$\begin{aligned}
 \|u\cdot \nabla u \|_{L_{10/7,1}(\R_+;L_{10/7}(\R^3))} &\leq C\|u\|_{L_\infty(\R_+;L_{3}(\R^3))} 
 \|\nabla u \|_{L_{10/7,1}(\R_+;L_{30/11}(\R^3))}\\ &\leq C\|u\|_{L_\infty(\R_+;L_{3}(\R^3))}
 \|\nabla^2u\|_{L_{10/7,1}(\R_+;L_{10/7}\R^3))}.\end{aligned}$$
Hence, thanks to~\eqref{eq:L31},
$$  \|u\cdot \nabla u \|_{L_{10/7,1}(\R_+;L_{10/7}(\R^3))} \leq C\|u_0\|_{L_2(\R^3)}^{2/3} \Xi^{1/3}\Psi.$$
Therefore, using~\eqref{t9}  and reverting to~\eqref{t10} yields
 $$ \Psi\leq C\bigl(\Psi_0 +  \|u_0\|_{L_2(\R^3)}^{2/3} \Xi_0^{1/3}\Psi\bigr),$$
 whence, thanks to the smallness condition in~\eqref{t9} (changing $C$ if need be),
 \begin{equation}\label{eq:PPsi}
 \Psi \leq 2C\Psi_0.\end{equation}
 Let us emphasize that, since $u_0$ is in $\dot B^{3/5}_{10/7,1}(\R^3)\cap\dot B^{6/5}_{5/2,1}(\R^3),$ it also belongs
 to all intermediate spaces, and in particular to $\dot B^{4/5}_{5/3,1}(\R^3)$ 
 with estimate $ \|u_0\|_{\dot B^{4/5}_{5/3,1}(\R^3)}\lesssim \|u_0\|_{\dot B^{3/5}_{10/7,1}(\R^3)}^{2/3}  \|u_0\|_{\dot B^{6/5}_{5/2,1}(\R^3)}^{1/3}.$
 Hence, mimicking the proof of \eqref{eq:PPsi}, we discover that, up to an irrelevant  change of $C,$ we have
 \begin{equation}\label{eq:WWW}
 \|u\|_{\dot W^{2,1}_{5/3,(5/3,1)}} \leq C \|u_0\|_{\dot B^{3/5}_{10/7,1}(\R^3)}^{2/3}  \|u_0\|_{\dot B^{6/5}_{5/2,1}(\R^3)}^{1/3}.
 \end{equation}

 \subsubsection*{Step 4. Time weight}

In order to get eventually  the desired control on $\nabla u$ in $L_1(\R_+;L_\infty(\R^3))$
that is needed to ensure~\eqref{eq:smallrho} provided we have~\eqref{eq:smallness} for $\rho_0,$ 
and, later on, uniqueness, 
we  mimic the sharp approach of the two dimensional case, 
considering the momentum equation in the following form
\begin{equation}\label{t11}
 (tu)_t - \Delta (tu) =(1-\rho)(tu)_t
 -t \rho u \cdot \nabla u + \rho u \mbox{ \ \ in \ }
 \R_+ \times \R^3.
\end{equation}
We observe that using Proposition \ref{p:MR} with $m=s=10/3$ delivers
\begin{equation}\label{t10a}\|u\|_{L_{10/3,1}(\R_+;L_{10/3}(\R^3))} \lesssim \Psi.
\end{equation}
Since the term $\rho u$ appears in the right-hand side of \eqref{t11}, it is natural 
  to look for a control  of $tu$ in the space
$\dot W^{2,1}_{10/3,(10/3,1)}(\R^3 \times \R_+).$
\smallbreak
Let $\Pi:= \sup_{t\geq0} \|tu\|_{\dot B^{7/5}_{10/3,1}} +  \| \nabla^2 (tu), (tu)_t\|_{L_{10/3,1}(\R_+;L_{10/3}(\R^3))}.$
 Proposition~\ref{p:MR}  and~\eqref{eq:smallrho} give us 
 $$\Pi\lesssim
 \|t u\cdot \nabla u \|_{L_{10/3,1}(\R_+;L_{10/3}(\R^3))} + \|u\|_{L_{10/3,1}(\R_+; L_{10/3}(\R^3))}.$$
In order to estimate the nonlinear term, we  first use  H\"older inequality to get:
\begin{equation}\label{t13}
 \| u\cdot \nabla t u \|_{L_{10/3,1}(\R_+ ; L_{10/3}(\R^3))} \leq C
 \|u\|_{L_{5,1}(\R_+ ; L_5(\R^3))} \|t\nabla u\|_{L_{10}(\R_+\times \R^3)}.
\end{equation}
{}From Gagliardo-Nirenberg inequality and Sobolev embedding, we know that 
$$\|\nabla z\|_{L_{10}(\R^3)}\lesssim \|\nabla^2z\|_{L_{10/3}(\R^3)}^{1/3}\|\nabla z\|_{\dot W^{2/5}_{10/3}(\R^3)}^{2/3}.$$
Therefore, using H\"older inequality, 
$$\|\nabla z\|_{L_{10}(\R_+\times \R^3)}\lesssim \|\nabla^2z\|_{L_{10/3}(\R_+\times\R^3)}^{1/3}\|z\|_{L_\infty(\R_+;\dot W^{7/5}_{10/3}(\R^3))}^{2/3}.$$
Consequently, \begin{equation}\label{t15}
 \|t\nabla u\|_{L_{10}(\R_+\times \R^3)} \lesssim \Pi.
\end{equation}
In order  to bound  $u$ in $L_{5,1}(\R_+;L_5(\R^3)),$
we use Proposition \ref{p:MR} with $p=q=5/3,$ and $m=s=5,$ and Inequality \eqref{eq:WWW}.  This gives
$$\|u\|_{L_{5,1}(\R_+;L_5(\R^3))}  \lesssim
\|u_0\|_{\dot B^{3/5}_{10/7,1}(\R^3)}^{2/3}  \|u_0\|_{\dot B^{6/5}_{5/2,1}(\R^3)}^{1/3}.$$
Putting together with~\eqref{t15} and reverting to~\eqref{t13}, we end up with 
$$
 \| u\cdot \nabla t u \|_{L_{10/3,1}(\R_+ \times \R^3)} \lesssim  \Psi_0^{2/3} \,\Xi_0^{1/3}\,\Pi.
 $$
Therefore, using also~\eqref{t10a} and \eqref{eq:PPsi}, we get the following inequality for $\Pi$:
$$\Pi \leq C\bigl(\Psi_0+ \Psi_0^{2/3} \,\Xi_0^{1/3}\,\Pi\bigr)\cdotp$$
Consequently, assuming 
$C \Psi_0^{2/3} \,\Xi_0^{1/3}\leq 1/2,$
a condition that implies~\eqref{t9} (up to a change of the constant maybe), 
we obtain 
\begin{equation}\label{t16} \Pi\leq 2 C\Psi_0.\end{equation}

 \subsubsection*{Step 5. Bounding $\nabla u$}  
It is now easy to get the desired control on $\nabla u$: we start from 
the following combination of the Gagliardo-Nirenberg inequality with Sobolev embedding: 
\begin{equation}\label{t17}
 \|\nabla u \|_{L_\infty(\R^3)} \leq C\|\nabla u\|^{2/3}_{\dot W^{1}_{10/3}(\R^3)} \|\nabla u\|^{1/3}_{
 \dot W^{1}_{5/2}(\R^3)},
\end{equation}
which implies 
\begin{equation}
 \int_0^\infty \|\nabla u\|_{L_\infty(\R^3)} dt\leq C\int_0^{\infty} t^{-2/3}
 \|t \nabla u\|_{\dot W^{1}_{10/3}(\R^3)}^{2/3} \|\nabla u \|^{1/3}_{\dot W^1_{5/2}(\R^3)}\,dt.
\end{equation}
Using H\"older inequality~\eqref{eq:holder}  with respect to time in Lorentz spaces, we find that
$$  \int_0^\infty \|\nabla u\|_{L_\infty(\R^3)}\, dt\leq C\|t^{-2/3}\|_{L_{3/2,\infty}(\R_+)}
  \|t\nabla u\|^{2/3}_{L_{10/3,1}(\R_+;\dot W^{1}_{10/3}(\R^3))} \|\nabla u \|^{1/3}_{L_{5/2,1}(\R_+;\dot W^1_{5/2}(\R^3))}.$$
As the right-hand side is bounded, owing to~\eqref{t9} and~\eqref{t16}, 
one may conclude that 
\begin{equation}\label{t18}
 \int_0^\infty \|\nabla u\|_{L_\infty(\R^3)}\, dt\leq C \Psi_0^{2/3}\, \Xi_0^{1/3}\ll 1.
\end{equation}
Therefore, arguing on the mass equation exactly as in the 2D case, 
one can justify~\eqref{eq:smallrho}, and thus all the previous steps
provided~\eqref{eq:smallness} is satisfied.

\subsubsection*{Step 6. Uniqueness}   Arguing as  in the previous section and 
knowing~\eqref{t18} (so as to put our system in Lagrangian coordinates),
 it suffices to establish  the additional 
property that $t^{1/2}\nabla u$ is in $L_2(0,T;L_\infty(\R^3)).$  Now, one may write, 
owing to~\eqref{t17}, that 
$$\begin{aligned}
\int_0^\infty t\|\nabla u\|_{L_\infty(\R^3)}^2\,dt
&\lesssim  \int_0^\infty t^{-1/3} \|t \nabla u\|_{\dot W^{1}_{10/3}(\R^3)}^{4/3} \|\nabla u \|^{2/3}_{\dot W^1_{5/2} (\R^3)}\,dt\\
&\lesssim \|t^{-1/3}\|_{L_{3,\infty}(\R_+)}
\bigl\|\|t \nabla u\|_{\dot W^{1}_{10/3}(\R^3)}^{4/3} \bigr\|_{L_{5/2,1}(\R_+)} 
\bigl\| \|\nabla u \|^{2/3}_{\dot W^1_{5/2} (\R^3)}\bigr\|_{L_{15/4}(\R_+)}\\
&\lesssim\|t \nabla^2 u\|_{L_{10/3,1}(\R_+; L_{10/3}(\R^3))}^{4/3}  
 \|\nabla^2 u \|^{2/3}_{L_{5/2}(\R_+\times\R^3)}.
\end{aligned}$$
Because $t\nabla^2 u$ is in $L_{10/3,1}(\R_+; L_{10/3}(\R^3))$
and $\nabla^2 u$ is in $L_{5/2}(\R_+\times\R^3),$ the right-hand
side is indeed bounded. This completes the proof of uniqueness. 

\subsubsection*{Step 7. Existence} Here we sketch the  proof of the existence 
of a global solution under our assumptions on the data. 
The overall strategy is essentially the same in dimensions $d=2$ and $d=3.$

As a first, we truncate $\rho_0$ and smooth out $u_0$ to meet the conditions
of  the local-in-time existence theorem of~\cite{DFP}: 
we  fix   a sequence $(\rho_0^n,u_0^n)_{n\in\N}$ that converges weakly to $(\rho_0,u_0)$
and satisfy the conditions therein. 
Let  $(\rho^n,u^n)_{n\in\N}$ be the corresponding sequence of maximal solutions, defined on $[0,T_n)\times\R^d$ and belonging for all $T<T_n$ to   the classical maximal regularity space
$$\dot W^{2,1}_{p,r}(T):= \bigl\{z\in \cC([0,T];\dot B^{2-2/r}_{p,r}(\R^d))\::\:
\d_tz,\nabla^2_xz\in L_r(0,T;L_p(\R^d))\bigr\},$$ 
with e.g.  $p=2d$ and $r=7/6.$
 \smallbreak
 It is shown in~\cite{DFP}  that those solutions  satisfy the energy balance
and~\eqref{eq:smallrho}.
Since the   computations of the previous step  just follow from the properties
of the heat flow and  basic functional analysis, each  $(\rho^n,u^n)$
satisfies the  estimates therein. 
In particular, $\|\nabla u^n\|_{L_1(0,T_n;L_\infty(\R^d))}$ is uniformly bounded like in~\eqref{t18}, 
which ensures control of~\eqref{eq:smallrho}.
Now, applying the standard maximal regularity estimates to\footnote{Of course
$\mu'$ is put to $0$ if one wants to prove the  existence part of Theorem~\ref{thm:NS3D}.}
$$
\d_tu^n-\Delta u^n-\mu'\Delta u^n=(1-\rho^n)\partial_tu^n+\rho^n u^n\cdot\nabla u^n,$$
one gets for all $T<T_n,$  
$$\|u^n\|_{\dot W^{2,1}_{p,r}(T)}\lesssim \|u_0^n\|_{\dot B^{2-2/r}_{p,r}(\R^d)}  + \|u^n\cdot\nabla u^n\|_{L_r(0,T;L_p(\R^d))},$$
whence 
\begin{equation}\label{eq:haut} \|u^n\|_{\dot W^{2,1}_{p,r}(T)}^r\lesssim \|u_0^n\|_{\dot B^{2-2/r}_{p,r}(\R^d)} ^r + 
\int_0^T \|u^n\|_{L_\infty(\R^d)}^r \|\nabla u^n\|_{L_p(\R^d)}^r\,dt.\end{equation}
Gagliardo-Nirenberg inequality reveals that
$$\|\nabla u^n\|_{L_p(\R^d)}\lesssim \|u^n\|_{\dot B^{2-2/r}_{p,r}(\R^d)}^{r/2}\|\nabla^2 u^n\|_{L_p(\R^d)}^{1-r/2}.$$
Therefore, plugging that inequality in \eqref{eq:haut} then using Young inequality, we discover that
for all $T<T_n$ and all $\varepsilon>0,$
$$\|u^n\|_{\dot W^{2,1}_{p,r}(T)}^r\leq C\|u_0^n\|_{\dot B^{2-2/r}_{p,r}(\R^d)}^r + 
C_\varepsilon \int_0^T \|u^n\|_{L_\infty(\R^d)}^2\|u^n\|_{\dot B^{2-2/r}_{p,r}(\R^d)}^r\,dt
+\varepsilon\int_0^T \|\nabla^2 u^n\|_{L_p(\R^d)}^r\,dt.$$
Then, taking $\varepsilon$ small enough and using Gronwall inequality, we end up with
$$
\|u^n\|_{\dot W^{2,1}_{p,r}(T)}^r\leq C\|u_0^n\|_{\dot B^{2-2/r}_{p,r}(\R^d)}^r
\exp\Bigl\{C\int_0^T\|u^n\|_{L_\infty(\R^d)}^2\,dt\Bigr\}\cdotp$$
Now, in the 2D case, we observe that 
$\|u^n\|_{L_\infty(\R^2)}\lesssim \|u^n\|_{L_2(\R^2)}^{1/3} \|\nabla^2u^n\|_{L_{4/3}(\R^2)}^{2/3},$
and thus, by H\"older inequality, 
$$\|u^n\|_{L_2(0,T_n;L_\infty(\R^2))}\lesssim \|u^n\|_{L_\infty(0,T_n;L_2(\R^2))}^{1/3} 
\|\nabla^2u^n\|_{L_{4/3}((0,T_n)\times\R^2)}^{2/3}\lesssim \|u_0\|_{\dot B^{1/2}_{4/3,1}(\R^2)},$$
and one can thus  bound   $u^n$ in $\dot W^{2,1}_{p,r}(T_n)$ independently  of $T_n.$ 
\medbreak
In the framework of Theorem~\ref{thm:NS3D}, we  use the following 
Gagliardo-Nirenberg inequality:
$$
\|u^n\|_{L_\infty(\R^3)}\lesssim  \|u^n\|_{\dot W^2_{5/2}(\R^3)}^{1/9}\|u^n\|_{\dot W^2_{10/7}(\R^3)}^{8/9}$$
that implies 
$$\|u^n\|_{L_{3/2}(\R_+;L_\infty(\R^3))} \lesssim \|u^n\|_{L_{5/2}(\R_+;\dot W^2_{5/2}(\R^3))}^{1/9}
\|u^n\|_{L_{10/7}(\R_+;\dot W^2_{10/7}(\R^3))}^{8/9},$$
then the fact that $L_2(\R_+;L_\infty(\R^3))\subset
L_{3/2}(\R_+;L_\infty(\R^3))\cap L_\infty(\R_+;L_\infty(\R^3))$
and that 
\newline $\dot B^{6/5}_{5/2,1}(\R^3)\hookrightarrow L_\infty(\R^3),$ 
to get the desired control  of $ \|u^n\|_{L_2(0,T;L_\infty(\R^3))}$ in terms
of $u_0$ only.   
\medbreak
In short, in both cases, one can  bound   $u^n$ in $\dot W^{2,1}_{p,r}(T_n)$ independently  of $T_n.$ 
Then, applying standard continuation arguments allows to prove that $(\rho^n,u^n)$ is actually global, 
and may be bounded in terms of the original data $(\rho_0,u_0)$ 
in the spaces of our main theorems, independently of $n.$

From this stage, passing to the limit \emph{in the slightly larger} (but reflexive) space 
$$\dot W^{2,1}_{5/2}(\R^3\times\R_+)\cap \dot W^{2,1}_{10/7}(\R^3\times\R_+)\quad 
\bigl(\hbox{or in }\  \dot W^{2,1}_{4/3}(\R^2\times\R_+)\bigr)$$ for the velocity 
can be done  as in~\cite{DFP} 
  (passing to the limit directly  in the nonreflexive spaces  $\dot W^{2,1}_{p,(p,1)}(\R^d\times\R_+)$
would require more care). 
The mass conservation equation may be handled according 
to   Di Perna and Lions' theory~\cite{DPL} (see details in~\cite{DFP})
and the momentum equation does not present any difficulty compared to works 
on global weak solutions, since a lot of regularity is available on the velocity
and there is no pressure term. 

Next, once  we know that $(\rho,u)$ is a solution, 
one can recover all the additional regularity, that
 are just based on `linear' properties
like interpolation or parabolic maximal regularity.\qed


\begin{appendix}
\section{}

Here we prove Proposition \ref{p:MR}.

 Performing the usual rescaling reduces the proof to $\mu=1.$
Now, the key element is the following interpolation relation proved in \cite[Th2:1.18.6]{Triebel}:  
$$ \left( L_{q_0,r_0}(X;A);L_{q_1,r_1}(X;A)\right)_{\theta,r}=L_{q,r}(X;A) \with\frac{1}{q}=\frac{1-\theta}{q_0}+\frac{\theta}{q_1}\andf \theta\in(0,1).$$
Taking $X=\R_+$ and $A=L_p(\R^d)$ thus leads to 
\begin{equation}\label{INTER}
 \left(L_{q_0}(\R_+;L_p(\R^d));L_{q_1}(\R_+;L_p(\R^d))\right)_{\theta,r}= L_{q,r}(\R_+;L_p(\R^d)).
\end{equation}
Now, based on the classical results for the heat equation,  one has the following maximal regularity estimates
for all $\alpha\in(1,\infty)$  and $1<p<\infty$: 
\begin{multline}\label{lplq}
 \|u\|_{L_\infty(\R_+;\dot B^{2-2/\alpha}_{p,\alpha}(\R^d))} + 
 \|u_t,\nabla^2 u\|_{L_{\alpha}(\R_+;L_{p}(\R^d))} \leq 
 C\left(  \|u_0\|_{\dot B^{2-2/\alpha}_{p,\alpha}(\R^d)} +\|f\|_{L_{\alpha}(\R_+;L_{p}(\R^d))}\right)\cdotp
\end{multline}
Let us take $\alpha=q_0,q_1,$ with $1<q_0<q<q_1<\infty$ such that $2/q=(1/q_0+1/q_1).$  Then,  
 the interpolation relation \eqref{INTER} ensures that 
$$
 \left(L_{q_0}(\R_+;L_{p}(\R^d));L_{q_1}(\R_+;L_{p}(\R^d))\right)_{1/2,r}=L_{q,r}(\R_+;L_{p}(\R^d))
$$
while the properties of interpolation for Besov spaces give us: 
$$
 \left(\dot B^{2-2/q_0}_{p,q_0}(\R^d);\dot B^{2-2/q_1}_{p,q_1}(\R^d)\right)_{1/2,r}=\dot B^{2-2/q}_{p,r}(\R^d).
$$
Hence, putting together  the above two relations with \eqref{lplq}  yields all the terms of \eqref{eq:maxreg1}, except for 
the norm in $L_{s,r}(\R_+;L_m(\R^d)).$ 
\medbreak
To achieve it, we observe that, provided $2-2/\alpha < d/p$,  Property \eqref{lplq} may be reformulated in the following terms: 
 $$\frac d{dt}-\Delta \ \hbox{ is an isomorphism  from the space }\ \dot W^{2,1}_{p,\alpha}(\R^d\times\R_+)\ \hbox{ onto }\ L_\alpha(\R_+;L_p(\R^d)).$$
Consequently, for  all $q_0$ and $q_1$ as above, $\frac d{dt}-\Delta$ is an isomorphism  from   
 $$ \bigl(\dot W^{2,1}_{p,q_0}(\R^d\times\R_+); \dot W^{2,1}_{p,q_1}(\R^d\times\R_+)\bigr)_{\frac12,r}
 \ \hbox{ onto }\ \bigl(L_{q_0}(\R_+;L_p(\R^d)); L_{q_1}(\R_+;L_p(\R^d))\bigr)_{\frac12,r}.$$
The latter space is $L_{q,r}(\R_+;L_{p}(\R^d)),$ and what we proved just above amounts to saying that 
$\frac d{dt}-\Delta$ is an isomorphism  from $\dot{W}^{2,1}_{p,(q,r)}(\R^d\times\R_+)$ to  $L_{q,r}(\R_+;L_{p}(\R^d)).$
Hence, we have 
\begin{equation}\label{eq:interpoW}
 \bigl(\dot W^{2,1}_{p,q_0}(\R^d\times\R_+); \dot W^{2,1}_{p,q_1}(\R^d\times\R_+)\bigr)_{\frac12,r} = \dot{W}^{2,1}_{p,(q,r)}(\R^d\times\R_+).\end{equation}
The end of the proof relies on the  mixed derivative theorem which ensures for each $\alpha\in(0,1)$ and $i=0,1,$ that
$$\dot W^{2,1}_{p,q_i}(\R^d\times\R_+) \hookrightarrow \dot W^\alpha_{q_i}(\R_+;\dot W^{2-2\alpha}_{p}(\R^d))$$
and on the following Sobolev embedding: 
$$\dot W^\alpha_{q_i}(\R_+;\dot W^{2-2\alpha}_{p}(\R^d))\hookrightarrow L_{s_i}(\R_+;L_m(\R^d))
\with \frac dm=\frac dp+2\alpha-2\andf \frac 1{s_i}=\frac 1{q_i}-\alpha.$$
Let us choose $\alpha=\frac 1q-\frac1s$ so that 
$$\frac12\biggl(\frac1{s_0}+\frac1{s_1}\biggr)=\frac12\biggl(\frac1{q_0}+\frac1{q_1}\biggr)-\alpha=\frac1s\cdotp$$
Since $(L_{s_0}(\R_+;L_m(\R^d)); L_{s_1}(\R_+;L_m(\R^d)))_{\frac12,r}=L_{s,r}(\R_+;L_m(\R^d)),$ 
this completes the proof of Inequality \eqref{eq:maxreg1}. 
\qed

\end{appendix}



\begin{thebibliography}{99}


\bibitem{Amann} H. Amann: {\em  Linear and quasilinear parabolic problems. Vol. I. Abstract linear theory.} Monographs in Mathematics, {\bf 89}, Birkh\"auser Boston, Inc., Boston, MA, 1995.



\bibitem{BCD} H. Bahouri, J.-Y. Chemin and  R. Danchin: {\it Fourier Analysis and Nonlinear Partial Differential Equations,} Grundlehren der mathematischen Wissenschaften, {\bf 343}, Springer (2011).

\bibitem{Be} C. Bennett, R. Sharpley: {\em Interpolation of operators.} Pure and Applied Mathematics, 129. Academic Press, Inc., Boston, MA, 1988.


 \bibitem{Ber}
 F. Berthelin: Existence and weak stability for a pressureless model with unilateral constraint.
 \emph{Math. Models Methods Appl. Sci.} {\bf 12} (2002), no. 2, 249--272. 

\bibitem{BDD}
F. Berthelin, P. Degond, M. Delitala and M. Rascle:
A model for the formation and evolution of traffic jams. 
\emph{Arch. Ration. Mech. Anal.} {\bf 187} (2008), no. 2, 188--220. 


\bibitem{BChZ}
D.  Bresch, Ch. Perrin and E. Zatorska: Singular limit of a Navier-Stokes system leading to a free/congested zones two-phase model. \emph{C. R. Math. Acad. Sci. Paris} {\bf 352} (2014), no. 9, 685--690.

\bibitem{CS}
 F. Cucker and  S. Smale: Emergent behavior in flocks. \emph{IEEE Trans. Automat. Control,} 
 {\bf 52} (2007), no. 5, 852--862.
 
\bibitem{Cwickel}
M. Cwickel: On $(L^{p_0} (A_0) , L^{p_1} (A_1))_{\theta , q}$, {\em Proc. Amer. Math. Soc.}, {\bf 44} (1974), 286--292.

\bibitem{D00} R. Danchin: Global existence in critical spaces for compressible
Navier-Stokes equations, {\em Inventiones Mathematicae}, {\bf 141}  (2000), no. 3,  579--614.

\bibitem{D-Fourier} R. Danchin: A Lagrangian approach for the compressible Navier-Stokes equations, {\em Annales de l'Institut Fourier}, {\bf 64} (2014), no. 2, 753--791.

\bibitem{DFP} R. Danchin, F. Fanelli and M. Paicu: A well-posedness result for viscous compressible fluids with only bounded density,  \emph{Analysis and PDEs}, 
{\bf 13} (2020), 275--316.

\bibitem{DM-adv}  R. Danchin and P.B. Mucha: 
    Compressible Navier-Stokes system : large solutions and incompressible limit, 
    \emph{Advances in Mathematics}, {\bf 320} (2017), 904--925.
    
   \bibitem{DM-1} R. Danchin and P.B. Mucha:  The incompressible Navier-Stokes equations in vacuum, {\em Communications
   on Pure and Applied Mathematics}, {\bf 52} (2019), 1351--1385.
  
  \bibitem{DM-TJM} R. Danchin and P.B. Mucha: 
  From compressible to incompressible inhomogeneous flows in the case of large data, 
  {\em Tunisian Journal of Mathematics}, {\bf 1} (2019), no. 1, 127--149.
  
  \bibitem{DM-Largenu} R. Danchin and P.B. Mucha: Compressible Navier-Stokes equations with ripped density (2019), arXiv:1903.09396. 
  
  \bibitem{DMPW}
R.  Danchin, P.B. Mucha, J. Peszek and B. Wr\'oblewski: Regular solutions to the fractional Euler alignment system in the Besov spaces framework.
{\em Math. Models Methods Appl. Sci.} {\bf 29} (2019), no. 1, 89--119.
  
  \bibitem{DBD} A. DeBenedictis and A. Das: {\em 
The General Theory of Relativity: A Mathematical Exposition.}
Springer, 2012.

\bibitem{Denk} R. Denk, M. Hieber, J. Pr\"uss:  R-boundedness, Fourier multipliers and problems of elliptic and parabolic type. Mem. Amer. Math. Soc. {\bf 166} (2003), no. 788.


\bibitem{DPL} R.  Di Perna and P.-L. Lions:
 Ordinary differential equations, transport theory and Sobolev spaces,
{\em Inventiones Mathematicae}, {\bf 98}  (1989) 511--547.

\bibitem{Kie}
T. Do, A. Kiselev, L. Ryzhik and C. Tan: Global regularity for the fractional Euler alignment system.
\emph{Arch. Ration. Mech. Anal.} {\bf 228} (2018) 1--37.

\bibitem{Fei}
E. Feireisl: {\em  Dynamics of Viscous Compressible Fluids,} Oxford Lecture Ser. Math. Appl. {\bf 26}, Oxford University Press, Oxford, 2004. 

\bibitem{FM} E. Feireisl, P.B. Mucha, A. Novotny, M. Pokorny: 
Time-periodic solutions to the full Navier-Stokes-Fourier system. 
{\em Arch. Ration. Mech. Anal. }{\bf 204}  (2012), no. 3, 745--786.

\bibitem{Giga} Y. Giga, H. Sohr: Abstract $L^p$ estimates for the Cauchy problem with applications to the Navier-Stokes equations in exterior domains. J. Funct. Anal. 102 (1991), no. 1, 72--94.

\bibitem{grafakos} L. Grafakos:
{\em Classical and Modern Fourier Analysis}, Prentice Hall, 2006.

\bibitem{HS}
M. Herty and V. Schleper:
Traffic flow with unobservant drivers.
\emph{ZAMM Z. Angew. Math. Mech.} {\bf 91} (2011), no. 10, 763--776. 

\bibitem{Hoff1} D. Hoff: Global solutions of the Navier-Stokes equations for multidimensional compressible flow with discontinuous initial data, {\em J. Differential Equations}, {\bf 120} (1995), no. 1, 215--254.

\bibitem{Hoff2} D. Hoff: Uniqueness of weak solutions of the Navier-Stokes equations of multidimensional, compressible flow,
{\em SIAM J. Math. Anal.}, {\bf 37} (2006), no. 6, 1742--1760.

\bibitem{Lions}
P.-L. Lions: \emph{Mathematical Topics in Fluid Mechanics. Vol. II: Compressible Models,} Oxford Lecture Ser. Math. Appl. {\bf 10}, The Clarendon Press, Oxford University Press, New York, 1998.

\bibitem{LSU}
O. Ladyzhenskaja, V. Solonnikov and N. Uraltseva: 
{\em Linear and quasilinear equations of parabolic type.}  Translations of Mathematical Monographs, {\bf 23} American Mathematical Society, Providence, R.I. 1967.

\bibitem{MN}
A. Matsumura and T. Nishida: {\em Initial boundary value problem for equations of motion of compressible viscous and heat conductive fluids}, Commun. Math. Phys. {\bf 89} (1983) 445--464. 

\bibitem{MM} P. Minakowski, P.B. Mucha, J. Peszek and  E. Zatorska: {\em  Singular Cucker-Smale dynamics.} Active particles, Vol. 2., 201--243, Model. Simul. Sci. Eng. Technol., Birkh\"auser/Springer, Cham, 2019.

\bibitem{MT} S. Motsch and E. Tadmor: A new model for self-organized dynamics and its flocking behavior,
 \emph{J. Stat. Phys.}  {\bf 144} (2011), 923--947.

 \bibitem{Mu03} P.B. Mucha:  The Cauchy problem for the compressible Navier-Stokes equations in the $L_p$-framework, {\em Nonlinear Anal.,} {\bf 52}, (2003), no. 4,  1379--1392. 
 
\bibitem{OPA} K-K. Oh, M-Ch. Park and  H-S. Ahn: A survey of multi-agent formation control.
\emph{Automatica,} {\bf  53} (2015), 424--440.

\bibitem{Pia} T. Piasecki, Y. Shibata, E. Zatorska: On the maximal $L^p$-$L^q$ regularity of solutions to a general linear parabolic system. J. Differential Equations {\bf 268} (2020), no. 7, 3332--3369.
 
\bibitem{Laura}
L. Saint-Raymond: \emph{Hydrodynamic limits of the Boltzmann equation.}
Lecture Notes in Mathematics, {\bf 1971}.  Springer-Verlag, Berlin, 2009.
 
  \bibitem{ST} R. Shvydkoy and  E. Tadmor: Eulerian dynamics with a commutator forcing.
 \emph{Trans. Math. Appl.} {\bf 1} (2017), no. 1, 26 pp.
 
\bibitem{TT}
E.  Tadmor and  Ch. Tan: Critical thresholds in flocking hydrodynamics with non-local alignment.
\emph{Philos. Trans. R. Soc. Lond. Ser. A  Math. Phys. Eng. Sci.} {\bf 372} (2014), no. 2028, 20130401, 22 pp. 

\bibitem{Triebel} H. Triebel:  {\em Interpolation theory, function spaces,
differential operators.} North-Holland Mathematical Library, {\bf 18}.
North-Holland Publishing Co., Amsterdam-New York, 1978.


\end{thebibliography}
\end{document}